\documentclass[12pt]{article}
\usepackage{amsmath,amsthm,amsfonts,amssymb,epsfig}
\usepackage{xcolor}

\usepackage{graphicx}

\usepackage{subcaption}

\newcommand{\reals}{\mathbb{R}}

\newtheorem{stel}{Theorem}
\newtheorem{gevolg}{Corollary}
\newtheorem{lemma}{Lemma}
\theoremstyle{remark}

\begin{document}

\title{Unique steady states for population models in a heterogeneous environment}
%\date{}
\author{Patrick De Leenheer\footnote{Department of Mathematics, Oregon State University, email: deleenhp@oregonstate.edu}, 
Jane S. MacDonald\footnote{Department of Mathematics, Oregon State University, email: macdojan@oregonstate.edu}, and 
Swati Patel\footnote{Department of Mathematics, Oregon State University, email: patelswa@oregonstate.edu}}
\date{}
\maketitle

\abstract{We revisit a model proposed by Freedman etal in \cite{freedman} which describes the dynamics of a population diffusing in a patchy environment. From their work it is known that positive steady states exist for this model, 
but not whether they are unique. Here, we provide sufficient conditions guaranteeing that steady states are unique. These conditions are satisfied when the reaction rates are generalized logistic growth rates. Our proofs critically exploit Chicone's ideas in \cite{chicone}, which were used to establish that the period map associated to a continuous family of periodic solutions of certain planar Hamiltonian systems is monotone.}

\section{Introduction}

Biological populations often reside in heterogeneous environments. In extreme cases the environment is patchy, and the populations may adjust their movement patterns based on the patch in which they currently reside. Furthermore, the habitat quality of the patches may differ significantly, affecting the population's growth or decay rates depending on where it resides. Here we consider a simple case of a two patch environment, and a population governed by a reaction-diffusion equation, in which the diffusion constants and reaction rates depend on the patch. 
We shall consider a corresponding steady state problem that takes the following form:
\begin{eqnarray}
0&=&\begin{cases}d^-u_{xx}+f^-(u),\textrm{ if } x\in (-L^-,0)\\ d^+u_{xx}+f^+(u),\textrm{ if }x\in (0,L^+) \end{cases} \label{pde} \\
u(0-)&=&u(0+)\textrm{ and } d^-u_x(0-)=d^+u_x(0+) \label{interface} \\
u_x(-L^-)&=&u_x(L^+)=0 \label{neumann}
\end{eqnarray}
This problem is defined on the spatial interval $[-L^-,L^+]$, where $L^-$ and $L^+$ are fixed, positive lengths. Equation $(\ref{pde})$ describes the steady states of a population with density $u(x)$, governed by distinct reaction-diffusion equations. They are characterized by the constant, positive diffusion constant $d^-$ and the density-dependent reaction term $f^-(u)$ on the left part $(-L^-,0)$ of the interval, and similarly by $d^+$ and $f^+(u)$ on the right part $(0,L^+)$. Condition $(\ref{interface})$ expresses that the density and the flux are continuous at the interface $x=0$; here, the notation $u(0-)$ stands for $\lim_{x\to 0-}u(x)$, and similarly for $u(0+)$, and for the derivatives $u_x(0-)$ and $u_x(0+)$. 
The boundary conditions $(\ref{neumann})$ express that there is no population flux across the boundary of the interval $[-L^-,L^+]$; this   corresponds to a Neumann boundary condition. 

Next we discuss the assumptions on the reaction terms $f^-(u)$ and $f^+(u)$. They represent the net growth rates of the population in the left and right parts of the interval respectively. 
We assume that the $f^{\pm}$ satisfy the following conditions:
\begin{enumerate}
\item $f^{\pm}$ are $C^1$ on $[0,+\infty)$, $f^{\pm}(0)=0$, and $(f^{\pm})'(0)>0$.
 
\item $f^{\pm}$ are $C^2$ on $(0,+\infty)$.

\item There exist distinct $K^{\pm}>0$ such that $f^{\pm}(u)\begin{cases}>0,\textrm{ if }u\in(0,K^{\pm})\\ 0,\textrm{ if }u=K^{\pm} \\ <0,\textrm{ if }u>K^{\pm} \end{cases}$.

\end{enumerate}
By reversing the orientation of the interval $[-L^-,L^+]$ if necessary, we can also assume without loss of generality that
\begin{equation}\label{hetero}
K^-<K^+,
\end{equation}
which reflects in part the heterogeneous nature of the environment.

For convenience, we introduce the following standing assumptions:
$$
{\bf SA} : \textrm{ The functions }f^{\pm} \textrm{ satisfy the above conditions 1,2 and 3, and } (\ref{hetero}) \textrm{ holds}.
$$

%Here used to be a Remark that proved that reaction rates satisfying ${\bf SA}$ could be represented alternatively by $f^{\pm}(u)=ug^{\pm}(u)$, for appropriate per capita growth rate functions $g^{\pm}(u)$. This "factorized" form is commonly used in the literature, see e.g. \cite{freedman,zaker}, which is why I decided to prove that our representation is in fact equivalent to this commonly used one, but in no way different. I have removed this proof and appended it at the very end of this document, so I can revisit if, if needed.

%\textcolor{red}{It would be nice if we could also allow for Gompertz growth rates, $f^{\pm}(u)=-r^{\pm}u\ln(u/K^{\pm})$. But these reaction rates are not differentiable at $u=0$. Note that they are continuous for $u=0$ however. Since all of the arguments involving Hamiltonian systems given below, are carried out for states where $u>0$, it is likely that we don't really need that the $f^{\pm}$ are $C^1$ on $[0,+\infty)$. If so, we could also allow for Gompertz growth rates here. We could check the details of this later. But on the other hand, condition ${\bf C2}$ always fails for Gompertz growth rates (I proved this in some handwritten notes), just like that condition fails for generalized logistic growth rates with $p^{\pm}\in(0,1)$.}

{\bf Example}: The main motivating example underlying the conditions imposed on the reaction rates $f^{\pm}$ is:
\begin{equation}\label{richards}
f^{\pm}(u)=r^{\pm}u\left(1-\left(\frac{u}{K^{\pm}} \right)^{p^{\pm}} \right),
\end{equation}
where the $r^{\pm}$, $K^{\pm}$ and $p^{\pm}$ are positive parameters. The $r^{\pm}$ are the maximal per capita growth rates, 
and the $K^{\pm}$ are the carrying capacities in the respective left and right intervals. These specific reaction rates are generalized logistic growth rate functions (known as Richards' functions in the literature) because they generalize the famous logistic growth rate functions (known as Verhulst's functions), obtained when the exponents $p^{\pm}$ are fixed to $1$.

Our main goal is to determine  if $(\ref{pde})-(\ref{neumann})$ has a unique, positive solution. More precisely, we say that 
$u(x)$, defined on $[-L^-,L^+]$, is a positive solution if $u(x)>0$ for all $x$ in $[-L^-,L^+]$, $C^2$ on $(-L^-,0)$ and $(0,L^+)$, $C^1$ on 
$[-L^-,0)$ and $(0,L^+]$, and such that $u(x)$ satisfies both the interface and boundary conditions $(\ref{interface})$ and $(\ref{neumann})$.

To address this problem, we shall first convert it into an equivalent ODE problem. We define two Hamiltonians $H^-(u,v)$ and $H^+(u,v)$, both defined on 
${\cal R}=\{(u,v)\in \reals^2 \,|\, u\geq 0\}$ :
\begin{equation}\label{ham}
H^{\pm}(u,v)=\frac{1}{2}v^2+F^{\pm}(u),\textrm{ where } F^{\pm}(u):=\frac{1}{d^{\pm}} \int_0^uf^{\pm}(s)ds.
\end{equation}
In the physics literature the functions $F^{\pm}(u)$ are referred to as the potential functions associated to the Hamiltonians $H^{\pm}$.  
The Hamiltonians $H^{\pm}$ give rise to two Hamiltonian systems on ${\cal R}$:
\begin{eqnarray}
u'&=& v, \label{ham-pm1}\\
v'&=& -\frac{1}{d^{\pm}}f^{\pm}(u), \label{ham-pm2}
\end{eqnarray}
where $'$ is shorthand notation for the derivative with respect to $x$. Level curves of the Hamiltonians $H^\pm$ correspond to orbits of $(\ref{ham-pm1})-(\ref{ham-pm2})$ in the $(u,v)$ phase plane  because $(H^\pm)'\equiv 0$. 
Note that $u(x)$, $x\in [-L^-,L^+]$, is a positive solution of $(\ref{pde})-(\ref{neumann})$ if and only if there exists a solution $\begin{pmatrix}u^-(x), v^-(x)\end{pmatrix}$, $x\in [-L^-,0]$, for the Hamiltonian system $(\ref{ham-pm1})-(\ref{ham-pm2})$ with Hamiltonian $H^-$ which starts at $x=-L^-$ at some $\begin{pmatrix}\alpha , 0 \end{pmatrix}$ where $\alpha>0$, and satisfies $u^-(x)>0$ for all $x\in [-L^-,0]$, and similarly a solution $\begin{pmatrix}u^+(x), v^+(x)\end{pmatrix}$, $x\in [0,L^+]$, for the Hamiltonian system $(\ref{ham-pm1})-(\ref{ham-pm2})$ with Hamiltonian $H^+$ which ends at $x=L^+$ at some $\begin{pmatrix}\beta , 0 \end{pmatrix}$ where $\beta>0$,  and satisfies $u^+(x)>0$ for all $x\in [0,L^+]$, and such that $u^-(0)=u^+(0)$, $d^-v^-(0)=d^+v^+(0)$, and 
$$
u(x)=\begin{cases}
u^-(x),\textrm{ if }x\in [-L^-,0],\\
u^+(x),\textrm{ if }x\in [0,L^+].
\end{cases}
$$

\begin{figure}
    \centering
    \includegraphics[scale = 0.3]{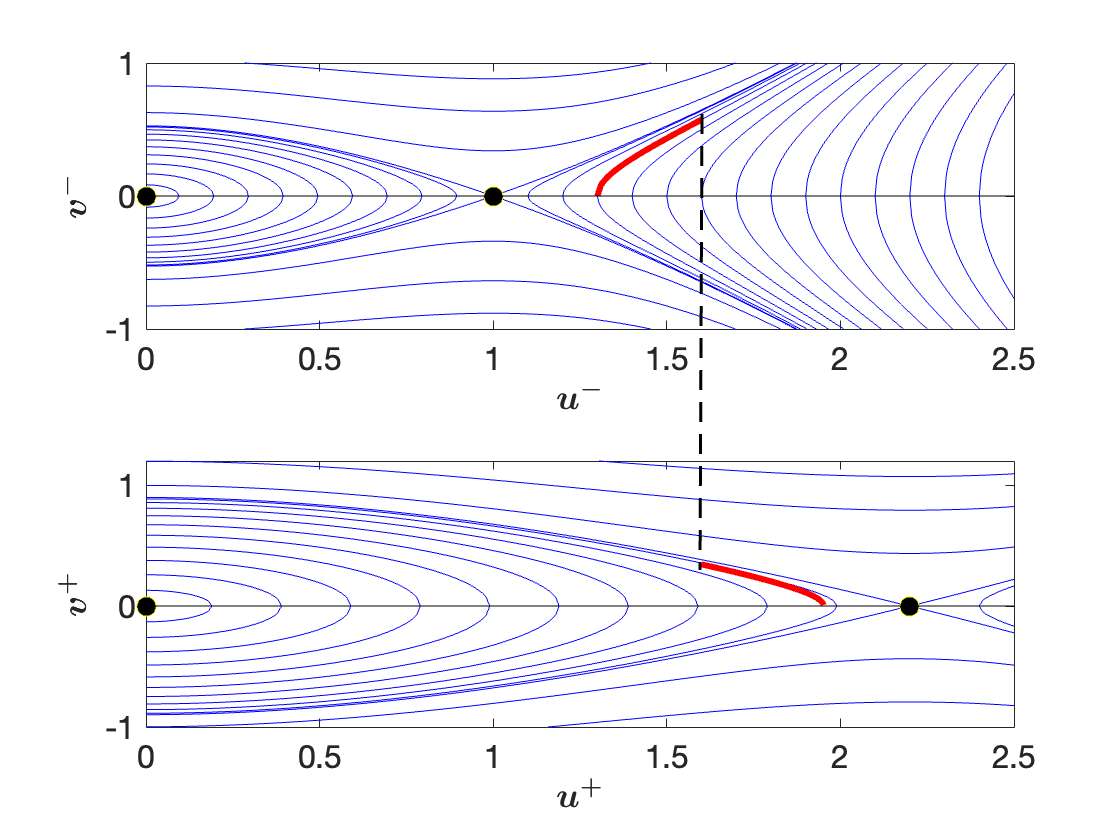}
    \caption{Piecewise orbits of the Hamiltonian systems $(\ref{ham-pm1})-(\ref{ham-pm2})$. Reaction rates $f^{\pm}(u)=r^{\pm}u(1-u/K^{\pm})$, with $r^{\pm}=1$, $K^-=1$ and $K^+=2.2$. Diffusion constants $d^-=1.2$ and $d^+=2$. Domain lengths $L^-=1.0349$, $L^+=1.1671$. Top and bottom orbits (in red) represent the solution $(u^-(x),v^-(x))$, $x\in [-L^-,0]$, of the $H^-$-system and the solution $(u^+(x),v^+(x))$, $x\in [0,L^+]$, of the $H^+$-system, respectively. The jump in the states at $x=0$ of these solutions occurs along the vertical dashed line (in black) and expresses the continuity of the density  $u^-(0)=u^+(0)$, and the flux $d^-v^-(0)=d^+v^+(0)$ at the interface. Other blue curves represent orbits of the Hamltonian systems, and the black dots their steady states.}
    \label{piecewise}
\end{figure}

Positive solutions of $(\ref{pde})-(\ref{neumann})$ can therefore be found by piecing together two specific solutions of the Hamiltonian systems $(\ref{ham-pm1})-(\ref{ham-pm2})$, as illustrated in Figure $\ref{piecewise}$. 
These consist of the $u$-coordinate $u^-(x)$, $x\in [-L^-,0]$, of a solution of the Hamiltonian system with Hamiltonian $H^-$ which starts at some point $(\alpha,0)$ on the positive $u$-axis at $x=-L^-$, and the $u$-coordinate $u^+(x)$, $x\in [0,L^+]$, of a solution of the Hamiltonian system with Hamiltonian $H^+$ which ends at some point $(\beta,0)$ at $x=L^+$, also on the positive $u$-axis; furthermore, both solutions must have 
positive $u$-coordinates for all $x$ in $[-L^-,L^+]$, and when evaluated at $x=0$, they must respect the interface condition $(\ref{interface})$. The strategy to find these solutions relies on the ``shooting method": Imagine two shooters who shoot arrows. One shooter is located at $(u,v)=(\alpha,0)$ and his arrow traces a forward solution of the Hamiltonian system with Hamiltonian $H^-$ that is $L^-$ units of ``time" underway; the other shooter is located at $(u,v)=(\beta,0)$ and his arrow traces a backward solution of the Hamiltonian system with Hamiltonian $H^+$ that is $L^+$ units of time underway. The problem amounts to finding a pair of positive numbers $\alpha, \beta$ which determine the locations of both shooters in such a way that their arrows end up in a point that respects the interface condition.

\section{Necessary conditions for positive solutions}

We begin by establishing some properties that any positive solution of $(\ref{pde})-(\ref{neumann})$ must have in Lemmas $\ref{nec1}$ and $\ref{nec2}$ below. These properties 
are most easily established by considering the phase portrait of the two Hamiltonian systems $(\ref{ham-pm1})-(\ref{ham-pm2})$, and although they 
were already obtained in \cite{freedman}, we include them here to make this paper self-contained. In so doing we also present a slightly different method of proof from those in \cite{freedman}.

\begin{lemma}\label{nec1}
Assume that ${\bf SA}$ holds. If $u(x)$, $x\in [-L^-,L^+]$, is a positive solution of $(\ref{pde})-(\ref{neumann})$, then 
\begin{equation*}
u(-L^-)>K^-,\textrm{ and }u(L^+)<K^+ .
\end{equation*}
\end{lemma}
\begin{proof}
The proof is by contradiction. Suppose that $u(-L^-)\leq K^-$. For all $x$ in $[-L^-,0]$, the positive solution $u(x)$ coincides with 
$u^-(x)$, for some solution $\begin{pmatrix}u^-(x), v^-(x) \end{pmatrix}$ of the Hamiltonian system  
$(\ref{ham-pm1})-(\ref{ham-pm2})$ with Hamiltonian $H^-$ starting at $x=-L^-$ at $\begin{pmatrix}\alpha , 0 \end{pmatrix}$ where $\alpha=u(-L^-)\leq K^-$.  
Then the solution $\begin{pmatrix}u^-(x), v^-(x) \end{pmatrix}$ remains in the set $B=\{(u,v)\in {\cal R}\,|\, 0<u\leq K^-,\; v\leq 0\}$ for all $x$ in $[-L^-,0]$, because ${\bf SA}$ holds. Therefore, 
$0<u^-(0)\leq K^-$ and $v^-(0)\leq 0$.

Now consider the solution $\begin{pmatrix}u^+(x), v^+(x) \end{pmatrix}$ of the Hamiltonian system  $(\ref{ham-pm1})-(\ref{ham-pm2})$ with Hamiltonian $H^+$ 
starting at $x=0$ in $\begin{pmatrix} u^-(0) , (d^-/d^+)v^-(0)\end{pmatrix}$. Then the positive solution $u(x)$ coincides with $u^+(x)$ for  all $x$ in $[0,L^+]$. But by ${\bf SA}$, $\begin{pmatrix}u^+(x), v^+(x) \end{pmatrix}$ remains in $B$ for  all $x$ in $[0,L^+]$ as well, and moreover, $v^+(x)<0$  for all $x\in (0,L^+]$. This implies in particular that $u_x(L^+)=v^+(L^+)<0$, contradicting that $u(x)$ satisfies the Neumann boundary condition $(\ref{neumann})$ at $x=L^+$. Therefore, $u(-L^-)>K^-$. 

An analogous proof shows that $u(L^+)<K^+$: Assuming that $u(L^+)\geq K^+$, backward  integration of the Hamiltonian system $(\ref{ham-pm1})-(\ref{ham-pm2})$ with Hamiltonian $H^+$ leads to 
a contraction in that the Neumann boundary condition $(\ref{neumann})$ at $x=-L^-$ cannot hold.

\end{proof}

\begin{lemma}\label{nec2}
Assume that ${\bf SA}$ holds. If $u(x)$, $x\in [-L^-,L^+]$, is a positive solution of $(\ref{pde})-(\ref{neumann})$, then $u_x(x)>0$ for all $x\in(-L^-,0) \cup (0,L^+)$, and $u_x(0-)$ and $u_x(0+)$ are positive, hence $u(x)$ is strictly increasing on $[-L^-,L^+]$. Furthermore, $u(x)\in (K^-,K^+)$ for all $x\in [-L^-,L^+]$.
\end{lemma}
\begin{proof}
Since $u(x)$, $x\in [-L^-,L^+]$ is a positive solution of $(\ref{pde})-(\ref{neumann})$, it coincides on $[-L^-,0]$ with $u^-(x)$, for some solution $\begin{pmatrix}u^-(x), v^-(x) \end{pmatrix}$ of the Hamiltonian system  
$(\ref{ham-pm1})-(\ref{ham-pm2})$ with Hamiltonian $H^-$ starting at $x=-L^-$ at $\begin{pmatrix}\alpha , 0 \end{pmatrix}$, where $\alpha:=u(-L^-)> K^-$ by Lemma $\ref{nec1}$.  Then the solution $\begin{pmatrix}u^-(x), v^-(x) \end{pmatrix}$ remains in the set ${\cal R}^-_\alpha:=\{(u,v)\in {\cal R}\,|\, u\geq \alpha,\; v\geq 0\}$ for all $x$ in $[-L^-,0]$, because ${\bf SA}$ holds, and moreover, that 
$v^-(x)>0$ for all $x\in(-L^-,0]$. Consequently, $u_x(x)=v^-(x)>0$ for all $x\in (-L^-,0)$ and $u_x(0-)=v^-(0)>0$. 

An entirely analogous proof shows that $u_x(x)>0$ for all $x\in (0,L^+)$ and that $u_x(0+)>0$.

The final statement of this Lemma follows from Lemma $\ref{nec1}$ and the fact that $u(x)$ is strictly increasing on $[-L^-,L^+]$.
\end{proof}

Lemmas $\ref{nec1}$ and $\ref{nec2}$ reveal that when assuming that only ${\bf SA}$ holds, any positive solution of $(\ref{pde})-(\ref{neumann})$ is necessarily strictly increasing, and takes values in the interval $(K^-,K^+)$. In the next section we focus 
on establishing conditions guaranteeing that a positive solution does indeed exist, and moreover, that this solution is unique.

\section{Sufficient conditions for unique positive solutions}

%\textcolor{red}{Here we should point out that any positive solution must necessarily lie in $[K^-,K^+]$. We should also still prove that the below map $u^-(0,\alpha)$ takes the value $K^+$ for some $\alpha^-$. This may require that we can guarantee that the solutions do not  blow up in finite time. A similar argument is then also needed to show there exists a $\beta^+$ for which $u^+(0,\beta^+)=K^-$.}

We shall study a class of parameterized solutions with parameter $\alpha$, of the Hamiltonian system $(\ref{ham-pm1})-(\ref{ham-pm2})$ having Hamiltonian $H^-$: For each $\alpha$ in $[K^-,K^+]$, we let $(u^-(x,\alpha),v^-(x,\alpha))$ denote the solution of this system starting at $x=-L^-$ in the initial condition $(\alpha,0)$, for all $x>-L^-$ for which this solution is defined.

We also impose a monotonicity condition on the reaction term $f^-$:
$$
{\bf M}^- : (f^-)'<0, \textrm{ on } [K^-,K^+]. 
$$
This condition aids in establishing that the maps $\alpha \to u^-(0,\alpha)$ and $\alpha \to v^-(0,\alpha)$ are strictly increasing in $\alpha$, as we show in the following Lemma.
\begin{lemma}\label{een}
Assume that ${\bf SA}$ and ${\bf M}^-$ hold. 

Then there exists some $\alpha^-\in (K^-,K^+)$ such that the map
$$
\alpha \to u^-(0,\alpha)
$$
is continuous and strictly increasing when  $\alpha\in [K^-,\alpha^-]$, and such that 
\newline $u^-(0,K^-)=K^-$ and $u^{-}(0,\alpha^-)=K^+$. In particular, this map takes $[K^-,\alpha^-]$ onto $[K^-,K^+]$. 
Furthermore, $u^-(0,\alpha) >K^+$ for $\alpha \in (\alpha^-,K^+]$ (provided that the solution $(u^-(x,\alpha),v^-(x,\alpha))$ did not blow up for some $x\in (-L^-,0)$).

In addition, the map
$$
\alpha \to v^-(0,\alpha)
$$
is continuous and strictly increasing on $[K^-,\alpha^-]$, and $v^-(0,K^-)=0$.

\end{lemma}
\begin{proof}
Consider the region ${\cal C}=\{(u,v)\in {\cal R} \,|\, K^-\leq u\leq K^+, \; v\geq 0\}$. Then the 
system $(\ref{ham-pm1})-(\ref{ham-pm2})$ with Hamiltonian $H^-$,  is  cooperative and irreducible \cite{smith} in ${\cal C}$, because ${\bf M}^-$ holds.  
The theory of monotone dynamical systems, see Theorem 4.1.1 in \cite{smith}, implies that solutions 
of system $(\ref{ham-pm1})-(\ref{ham-pm2})$ in ${\cal C}$ that have distinct, ordered initial conditions, are strongly ordered for all future times, as long as they remain in ${\cal C}$. This implies that the maps $\alpha \to u^-(0, \alpha)$ and $\alpha \to v^-(0,\alpha)$ are strictly increasing, provided that $(u^-(x,\alpha),v^-(x,\alpha))$ remains in ${\cal C}$ for all $x$ in $[-L^-,0]$. Since $(K^-,0)$ is a steady state 
of the Hamiltonian system, this steady state solution certainly remains in ${\cal C}$, and we have that $u^-(0,K^-)=K^-$ and $v^-(0,K^-)=0$. Continuity of the maps $\alpha \to u^-(0,\alpha)$ and $\alpha \to v^-(0,\alpha)$ --which follows from the basic qualitative theory of ODE's-- implies that for sufficiently small $\alpha>K^-$, the solution $(u^-(x,\alpha),v^-(x,\alpha))$ also remains in ${\cal C}$ for all $x \in [-L^-,0]$. On the other hand, by ${\bf SA}$, we have that $u^-(0,K^+)>K^+$, provided that the solution starting in $(0,K^+)$ at $x=-L^-$ did not blow up for some $x$ in $(-L^-,0)$. By the continuity of 
$u^-(0,\alpha)$ with respect to $\alpha$, it follows from the intermediate value theorem that  
there exists an $\alpha^-$ in $(K^-,K^+)$ such that $u^-(0,\alpha^-)=K^+$, and furthermore, that $v^-(0,\alpha^-)>0$, which follows from ${\bf SA}$. 
And since the map $\alpha \to u^-(0,\alpha)$ is strictly increasing on the interval $[K^-,\alpha^-]$, and $u^-(0,\alpha^-)=K^+$, we have that $u^-(0,\alpha)<K^+$ for all $\alpha$ in $[K^-,\alpha^-)$. 
Finally, note that $u^-(0,\alpha)>K^+$ for all $\alpha$ in $(\alpha^-,K^+]$ (provided the corresponding solution of the Hamiltonian system did not blow up). It follows that $\alpha^-$ is the unique value of $\alpha$ in $(K^-,K^+)$ for which $u^-(0,\alpha^-)=K^+$. 
\end{proof} 

Next we study a class of parameterized solutions having parameter $\beta$, of the Hamiltonian system $(\ref{ham-pm1})-(\ref{ham-pm2})$ with Hamiltonian $H^+$: For each $\beta$ in $[K^-,K^+]$, we let $(u^+(x,\beta),v^+(x,\beta))$ denote the (backward) solution of this system starting at $x=+L^+$ in the initial condition $(\beta,0)$, for all $x<+L^+$ for which this solution is defined.

We impose two conditions on the potential function $F^+$, which was defined in $(\ref{ham})$, associated to the Hamiltonian $H^+$. These conditions express the concavity and convexity respectively of two functions related to $F^+$:
\begin{eqnarray*}
{\bf C1}^+&:& (\sqrt{F^+})''\leq 0,\textrm{ on } (K^-,K^+). \\
{\bf C2}^+&:& \left(\frac{F^+}{((F^+)')^2} \right)'' \geq 0, \textrm{ on }(K^-,K^+).
\end{eqnarray*}

\begin{lemma} \label{twee}
Assume that ${\bf SA}$, ${\bf C1}^+$ and ${\bf C2}^+$ hold.
  
Then there exists some $\beta^+ \in (K^-,K^+)$ such that the map
$$
\beta \to u^+(0,\beta)
$$
is continuous and strictly increasing when $\beta \in [\beta^+,K^+]$, and such that 
\newline $u^+(0,\beta^+)=K^-$ and $u^{+}(0,K^+)=K^+$. In particular, this map takes $[\beta^+,K^+]$ onto $[K^-,K^+]$.
Furthermore, $u^+(0,\beta) <K^-$ for $\beta \in (K^-,\beta+]$ (provided that the solution $(u^+(x,\beta),v^+(x,\beta))$ remains in ${\cal R}$ for all $x\in [0,L^+]$).

In addition, the map
$$
\beta \to v^+(0,\beta)
$$
is continuous and strictly decreasing on $[\beta^+,K^+]$,  and $v^+(0,K^+)=0$.

\end{lemma}
\begin{proof}
The proof is similar to the proof of Lemma $\ref{een}$, but differs from it in one important sense: To show that the maps $\beta \to u^+(0,\beta)$ and $\beta \to v^+(0,\beta)$ are strictly increasing and strictly decreasing respectively on $[\beta^+,K^+]$, we can no longer exploit the theory of monotone dynamical systems. Indeed, the Hamiltonian system $(\ref{ham-pm1})-(\ref{ham-pm2})$ is not necessarily a monotone and irreducible system in the region ${\cal C}=\{(u,v)\in {\cal R} \,|\, K^-\leq u\leq K^+, \; v\geq 0\}$. However, Lemmas $\ref{monotone-period1}$ and $\ref{monotone-period2}$,  which are proved in the Appendix, enable us to overcome this problem and both Lemmas imply that  
the maps $\beta \to u^+(0,\beta)$ and $\beta \to v^+(0,\beta)$ are strictly increasing and decreasing, respectively. 
%That the maps $u^+(0,\beta)$ and $v^+(0,\beta)$ evaluate to the claimed values when $\beta=\beta^+$ and $\beta=K^+$, follows from basic properties of the (backward) flow of the Hamiltonian system $(\ref{ham-pm1})-(\ref{ham-pm2})$ having Hamiltonian $H^+$, which  include that 
%$(K^+,0)$ is a steady state of this system. 

\end{proof}
By combining Lemmas $\ref{een}$ and $\ref{twee}$ we obtain
\begin{gevolg} \label{cor} Assume that ${\bf SA}$, ${\bf M}^-$, ${\bf C1}^+$ and ${\bf C2}^+$ hold. 

Then there is a unique $(\alpha^*,\beta^*)$ in $(K^-,\alpha^-)\times (\beta^+,K^+)$ such that
\begin{eqnarray}
u^-(0,\alpha^*)&=&u^+(0,\beta^*) \label{den}\\
d^-v^-(0,\alpha^*)&=&d^+v^+(0,\beta^*) \label{flu}
\end{eqnarray}
\end{gevolg}

\begin{proof}
First, since the two maps $\alpha \to u^-(0,\alpha)$, defined for $\alpha\in [K^-,\alpha^-]$, and $\beta \to u^+(0,\beta)$, defined for $\beta \in [\beta^+,K^+]$, are continuous, strictly increasing and onto $[K^-,K^+]$ by Lemmas ${\ref{een}}$ and ${\ref{twee}}$, it follows that there is a unique, continuous and strictly increasing  
function $\alpha \to \beta(\alpha)$, defined for all $\alpha$ in $[K^-,\alpha^-]$, and onto the interval $[\beta^+,K^+]$, such that
$$
u^-(0,\alpha)=u^+(0,\beta(\alpha)).
$$
Next, consider the continuous map $\alpha \to d^+v^+(0,\beta(\alpha))-d^-v^-(0,\alpha)$, for $\alpha$ in $[K^-,\alpha^-]$. Lemmas $\ref{een}$ and $\ref{twee}$, and the fact that $\beta(\alpha)$ is strictly increasing, imply that this map is strictly decreasing. Furthermore, it takes a positive value when $\alpha=K^-$, and a negative value when $\alpha=\alpha^-$. Thus, by the intermediate value theorem, there is a unique $\alpha^*$ in $(K^-,\alpha^-)$ and a unique $\beta^*:=\beta(\alpha^*)$ in $(\beta^+,K^+)$ such that $(\ref{den})-(\ref{flu})$ hold.
\end{proof}

In summary, by combining Lemmas $\ref{nec1}$, $\ref{nec2}$, $\ref{een}$, $\ref{twee}$ and Corollary $\ref{cor}$, we have proved the first main result of this paper:
\begin{stel}\label{main}
Assume that ${\bf SA}$, ${\bf M}^-$, and ${\bf C1}^+$ and ${\bf C2}^+$ hold. 

Then $(\ref{pde})-(\ref{neumann})$ has a unique, positive solution.
\end{stel}

{\bf Application}: Below we identify a broad class of reaction terms $f^{\pm}(u)$, consisting of certain generalized logistic growth rate functions, for which all the conditions in Theorem $\ref{main}$ hold.
\begin{lemma} \label{application}
Assume that $f^{\pm}(u)$ is given by $(\ref{richards})$, and that $K^-<K^+$. 

Then clearly, {\bf SA} hold. Furthermore,
\begin{enumerate}
\item ${\bf M}^-$ holds,

\item ${\bf C1}^+$ holds for all $K^-$ in $(0,K^+)$, and

\item ${\bf C2}^+$ holds for all $K^-$ in $(0,K^+)$ if and only if $p^+\geq 1$.

\end{enumerate}
\end{lemma}

\begin{proof}
The first item in this Lemma is straightforward to check. 
The proof of the second and third items in this Lemma are also straightforward but somewhat lengthy exercises in differential calculus. We only highlight the main steps.

We begin by showing that ${\bf C1}^+$ holds for all $K^-$ in $(0,K^+)$. Using $(\ref{richards})$, this happens
if and only if
$$
(r^+)^2u^2Q((u/K^+)^{p^+})\leq 0, \textrm{ for all } u\in (0,K^+),
$$
where 
$$
Q(z):=\left(1-\frac{2}{p^++2}z \right)(1-(p^++1)z)- (1-z)^2.
$$

Thus, to verify that ${\bf C1}^+$ holds for all $K^-$ in $(0,K^+)$, it suffices to show that the quadratic polynomial $Q(z)$ is non-positive for all $z$ in $[0,1]$. After some algebra, one finds that $Q(z)$ factors as follows:
$$
Q(z)=\frac{p^+}{p^++2}z(z-(p^++1)).
$$
From this it is clear that $Q(z)\leq 0$ for all $z$ in $[0,1],$ because $p^+$ is positive. 

Next we show that ${\bf C2}^+$ holds for all $K^-$ in $(0,K^+)$  if $p^+\geq 1$. Using $(\ref{richards})$, 
$$
\frac{F^+(u)}{((F^+)'(u))^2}=:R\left(g(u)\right),\textrm{ where }R(z):=\frac{1}{2r^+}\frac{1-\frac{2}{p^++2}z}{(1-z)^2}\textrm{ and } g(u):=\left(\frac{u}{K^+}\right)^{p^+}.
$$
Observe that $g(u)\in (0,1)$ for all $u$ in $(0,K^+)$. This representation of 
\newline $F^+/((F^+)')^2$ implies that
$$
\left(\frac{F^+(u)}{((F^+)'(u))^2}\right)'' = R''(g(u))(g'(u))^2+R'(g(u))g''(u)
$$
Note that $g'(u)=(p^+/K^+)(u/K^+)^{p^+-1}$, and 
\newline $g''(u)=(p^+(p^+-1)/(K^+)^2)(u/K^+)^{p^+-2}$. 
Thus, if $p^+\geq 1$, then $g''(u)\geq 0$ for all $u$ in $(0,K^+)$. Therefore, to show that ${\bf C2}^+$ holds for all $K^-$ in $(0,K^+)$ when $p^+\geq 1$, it suffices to prove that $R'(z)\geq 0$ and $R''(z)\geq 0$, for all $z$ in $(0,1)$. 
Straightforward  calculations reveal that
$$
R'(z)=\frac{1}{r^+(p^++2)}\frac{-z+(p^++1)}{(1-z)^3},\textrm{ and }R''(z)=\frac{1}{r^+(p^++2)} \frac{-2z+(3p^++2)}{(1-z)^4}.
$$
Thus $R'(z)$ and $R''(z)$ are positive for all $z$ in $(0,1)$.

To conclude this proof, we will show that ${\bf C2}^+$ fails if $p^+\in (0,1)$. Suppressing some tedious algebra,
$$
\left(\frac{F^+(u)}{((F^+)'(u))^2}\right)''=:\frac{p^+}{r^+(K^+)^2(p^++2)(1-(u/K^+)^{p^+})^4}\left(\frac{u}{K^+} \right)^{p^+-2}P\left(\frac{u}{K^+}\right),
$$
$$
\textrm{ where }P(z):=p^+z((3p^+ +2)-2z)+(p^+ -1)((p^+ +1)-z)(1-z).
$$
To show that ${\bf C2}^+$ fails for  some $K^-$ in $(0,K^+)$ when $p^+\in (0,1)$, it suffices to prove that the quadratic polynomial $P(z)$ changes sign in the interval $[0,1]$. 
This sign change occurs because
$$
P(0)=(p^+)^2-1<0,\textrm{ and }P(1)=3(p^+)^2>0.
$$

\end{proof}
\section{An alternative sufficient condition}
In some cases, it may happen that $(f^-)'$ changes sign in the interval $[K^-,K^+]$. Then ${\bf M}^-$ is not satisfied and we cannot use Theorem $\ref{main}$ to establish the existence of a unique positive solution for system $(\ref{pde})-(\ref{neumann})$. Fortunately, in such cases, we can replace the monotonicity condition ${\bf M}^-$ by two conditions imposed on the potential function $F^-$ associated to the Hamiltonian $H^-$, which still imply the existence of a unique positive solution. These conditions are similar to ${\bf C1}^+$ and ${\bf C2}^+$. To introduce them, 
%As before, these conditions express the concavity and convexity respectively of two functions related to $F^-$. 
first recall that if ${\bf SA}$ holds, then $K^-<K^+$, and $(F^-)'=f^-<0$ on $(K^-,K^+]$. Therefore, the function 
\begin{equation}\label{new-potential}
G^-(u):=F^-(u)-F^-(K^+),
\end{equation}
is positive for all $u$ in $(K^-,K^+)$. We now introduce a pair of conditions for the function $G^-$:
\begin{eqnarray*}
{\bf C1}^-&:& (\sqrt{G^-})''\leq 0,\textrm{ on } (K^-,K^+). \\
{\bf C2}^-&:& \left(\frac{G^- }{((G^-)')^2} \right)'' \geq 0, \textrm{ on }(K^-,K^+).
\end{eqnarray*}

We claim that the conclusions of Lemma $\ref{een}$ remain valid when condition ${\bf M}^-$ is replaced by the pair of conditions ${\bf C1}^-$ and ${\bf C2}^-$. To see why, we revisit the proof of Lemma $\ref{een}$ and modify it following the arguments given in the proof of Lemma $\ref{twee}$. But in this case, to establish the key fact that 
the maps $\alpha \to u^-(0,\alpha)$ and $\alpha \to v^-(0,\alpha)$ are strictly increasing when $\alpha$ belongs to $[K^-,\alpha^-]$, we rely on Lemmas $\ref{monotone-period3}$ and $\ref{monotone-period4}$ (instead of Lemmas $\ref{monotone-period1}$ and $\ref{monotone-period2}$ as we did in the proof of Lemma $\ref{twee}$); Lemmas $\ref{monotone-period3}$ and $\ref{monotone-period4}$ are proved in the Appendix.

The second main result of this paper below then follows by combining Lemmas $\ref{nec1}$, $\ref{nec2}$, $\ref{een}$ (with condition ${\bf M}^-$ replaced by the pair ${\bf C1}^-$ and ${\bf C2}^-$), 
$\ref{twee}$ and Corollary $\ref{cor}$.

\begin{stel}\label{main2}
Assume that ${\bf SA}$, ${\bf C1}^-$, ${\bf C2}^-$, ${\bf C1}^+$ and ${\bf C2}^+$ hold. 

Then $(\ref{pde})-(\ref{neumann})$ has a unique, positive solution.
\end{stel}

\section{Conclusion}
We have considered a population model governing the density of a species in a two patch environment. The model was first investigated in \cite{freedman} where it was shown that positive steady states always exist. This raises the natural question of whether or not positive steady states are in fact unique. Uniqueness of steady states has been claimed more recently, see Theorem 1 in \cite{zaker}, but it appears that there is a gap in the proof. In this paper, we have focused on finding sufficient conditions which guarantee that the model has a unique, positive steady state. 
Positive steady states are unique, provided that certain maps associated to particular Hamiltonian systems are strictly monotone. We have proposed  
sufficient conditions on the reaction terms (or on their integrals) which guarantee that these maps are indeed strictly monotone. 
These sufficient conditions were inspired by Chicone's work \cite{chicone}, who used them to prove the strict monotonicity of a period map  encountered in the study of certain planar Hamiltonian systems. 
While our sufficient conditions are satisfied for a broad class of generalzed logisitic growth rate functions $(\ref{richards})$, namely when $p^{\pm}\geq 1$, they can fail when the exponents $p^{\pm}$ belong to the interval $(0,1)$, at least for certain combinations of the carrying capacities $K^-$ and $K^+$, see Lemma $\ref{application}$. It is unclear at this time whether the failure of our sufficient conditions in this case implies that the model could have multiple positive steady states, or whether the steady state is still unique.
\newpage

\section*{Appendix}
We consider the planar Hamiltonian system with Hamiltonian $H^+(u,v)=v^2/2+F^+(u)$. Since $F^+(u)=(1/d^+)\int_0^u f^+(s)ds$, the standing assumptions {\bf SA} imply that the potential $F^+(u)$ satisfies the following conditions:
\begin{enumerate}
\item $F^+$ is $C^2$ on $[0,+\infty)$, $F(0)=F'(0)=0$, and $F''(0)>0$.
\item $F^+$ is $C^3$ on $(0,+\infty)$.
\item There is a $K^+>0$ such that $(F^+)'(u)\begin{cases} >0,\textrm{ if } u\in(0,K^+),\\ =0,\textrm{ if }u=K^+,\\ <0,\textrm{ if }u>K^+. \end{cases}$
\end{enumerate}
%For notational convenience we denote the derivative $(F^+)'(u)$ by $f^+(u)$ in what follows.

%{\bf Example}:  We are particularly interested in the potential function
%\begin{equation}\label{potential}
%F^+(u)=ru^2\left(\frac{1}{2}-\frac{1}{p+2}\left(\frac{u}{K}\right)^p \right),
%\end{equation}
%where $r^+,K^+$ and $p^+$ are positive parameters. Then
%\begin{equation}\label{gen-logistic}
%f^+(u)=(F^+)'(u)=r^+u\left(1-\left(\frac{u}{K^+} \right)^{p^+} \right),
%\end{equation}
%is the generalized logistic growth rate (aka the Richards function). It generalizes the logistic growth rate (or Verhulst function) $^+ru(1-u/K^+)$, which appears in many models of biological systems.

Our focus will be on a certain solutions of
the Hamiltonian system 
\begin{eqnarray}
u'&=& v, \label{ham1}\\
v'&=& -\frac{1}{d^+}f^+(u), \label{ham2}
\end{eqnarray}
defined in the region ${\cal R}=\{(u,v)\in \reals^2\,| \, u\geq 0\}$. Since the Hamiltonian $H^+(u,v)$ is conserved along solutions of $(\ref{ham1})-(\ref{ham2})$ (indeed, $(H^+)'\equiv 0$), the orbits of the solutions of this system are contained in the level curves of the Hamiltonian. 

Fix any $u_0$ in $(K^-,K^+)$ and define $E_{u_0}:=F^+(u_0)$ and $E_{K^+}:=F^+(K^+)$. Since $F^+$ is strictly increasing in $(0,K^+)$, and $F^+(0)=0$, it follows that $0<E_{u_0}<E_{K^+}$. Furthermore, by the Implicit Function Theorem, there is a unique $C^3$ function 
$E\to \beta(E)$, defined for all $E$ in $(E_{u_0},E_{K^+})$, such that $F^+(\beta(E))=E$ and $\beta(E)\in (u_0,K^+)$.

\begin{figure}
    \centering
    \includegraphics[scale = 0.2]{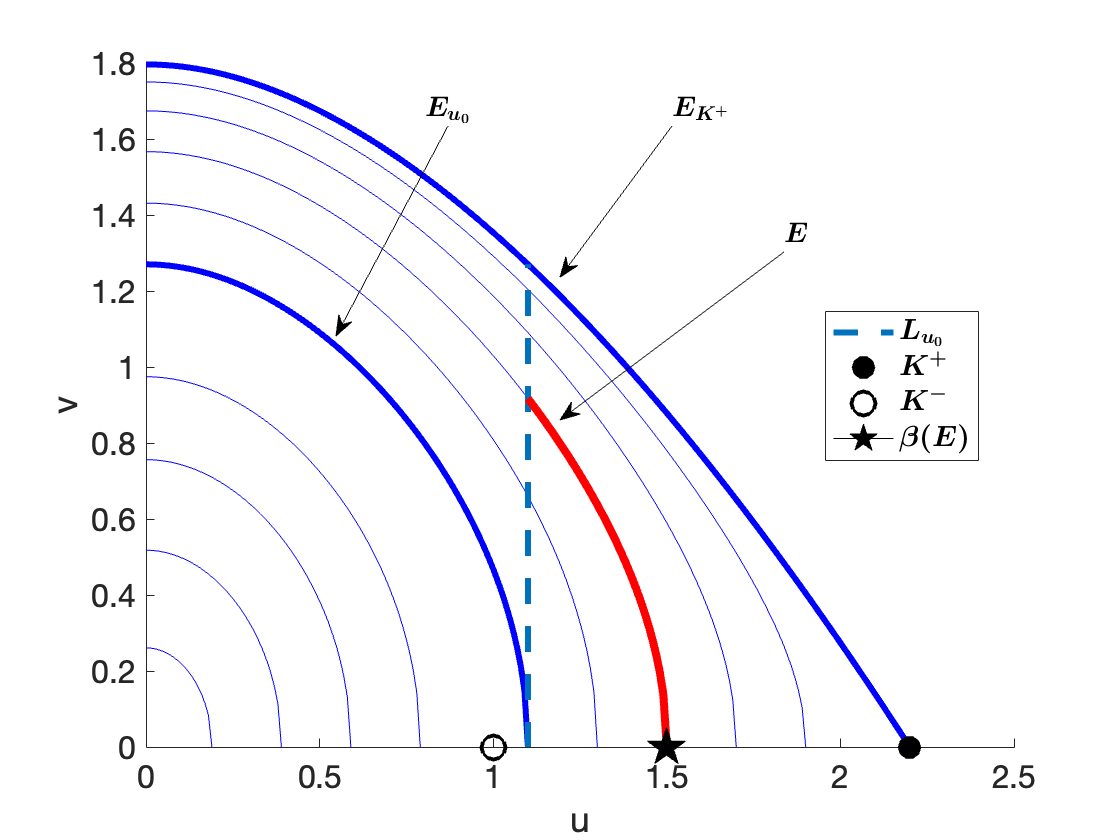}
    \caption{Orbits of the Hamiltonian system $(\ref{ham1})-(\ref{ham2})$. Reaction rate $f^+(u)=r^+u(1-u/K^+)$, with $r^+=1$, $K^+=2.2$ and diffusion constant $d^+=2$. Line segment $L_{u_0}$ for 
    $u_0=1.1$.}
    \label{transversal_u0}
\end{figure}

We shall consider a family of solutions of $(\ref{ham1})-(\ref{ham2})$ whose initial conditions belong to the line segment $L_{u_0}=\{(u_0,v)\,|\, 0<v<\sqrt{2(E_{K^+}-E_{u_0})}\}$, as illustrated in Figure $\ref{transversal_u0}$. 
The vector field of $(\ref{ham1})-(\ref{ham2})$ is transversal to $L_{u_0}$. 
We parameterize the solutions starting on $L_{u_0}$ by the (constant) value $E$ the Hamiltonian $H^+$ takes on them, and note that $E$ ranges over $(E_{u_0},E_{K^+})$ as the initial condition of the solution varies in $L_{u_0}$. 
Note that for each $E$ in $(E_{u_0},E_{K^+})$, the corresponding solution intersects the $u$-axis in the point $(\beta(E),0)$. 
For each $E$ in $(E_{u_0},E_{K^+})$, let $T^+_{u_0}(E)$ 
denote the time it takes for the solution to reach the point $(\beta(E),0)$. 
%The map $T(E)$ is (at least) $C^1$ on the interval $(E_{u_0},E_K)$ because the vector field of system $(\ref{ham1})-(\ref{ham2})$ is $C^2$ in ${\cal R}$. In fact, 
Then $T^+_{u_0}(E)$ can be represented by the following Lebesgue integral:
\begin{equation}\label{period}
T^+_{u_0}(E)=\int_{u_0}^{\beta(E)} \frac{du}{\sqrt{2(E-F^+(u))}},\textrm{ for all } E\in (E_{u_0},E_{K^+}).
\end{equation}
To see why $T^+_{u_0}(E)<+\infty$, for all $E$ in $(E_{u_0},E_{K^+})$, we first fix $E$. Since $(F^+)'$ is positive and continuous on $[u_0,\beta(E)]$, $(F^+)'$ achieves a positive minimum, say $m_E$, on this compact interval. It follows from the Mean Value Theorem that
$$
F^+(\beta(E))-F^+(u)\geq m_E(\beta(E)-u),\textrm{ for all }u\in [u_0,\beta(E)].
$$
Recalling that $F^+(\beta(E))=E$, it follows that
$$
T^+_{u_0}(E)\leq \frac{1}{\sqrt{2m_E}}\int_{u_0}^{\beta(E)}\frac{du}{\sqrt{\beta(E)-u}} <+\infty.
$$
Our main goal is to ensure that $T^+_{u_0}(E)$ is strictly increasing with respect to $E$. 
\begin{lemma} \label{monotone-period1}
Suppose that ${\bf SA}$, ${\bf C1}^+$ and ${\bf C2}^+$ hold. Fix any $u_0$ in $(K^-,K^+)$. 

Then $\frac{dT^+_{u_0}(E)}{dE}>0$, for all $E$ in $(E_{u_0},E_{K^+})$.
\end{lemma}
\begin{proof}
In order to determine the derivative of $T^+_{u_0}(E)$, we will follow Chicone's strategy \cite{chicone} by first converting 
the integral in $(\ref{period})$ by means of two appropriate substitutions.
To begin this process, we first define
\begin{equation}\label{h-function}
h(u)=(F^+(u))^{1/2},\textrm{ for all } u\in (0,K^+). \\
%\begin{cases}
%u\left(\frac{F(u)}{u^2} \right)^{1/2},\textrm{ if } u>0 \\
%0,\textrm{ if }u=0
%\end{cases}
\end{equation}
Then $h'(u)>0$ for all $u\in (0,K^+)$, and therefore the inverse $h^{-1}$ of $h$ exists. Also note that $h$ is $C^3$ on $(0,K^+)$, and so is $h^{-1}$ on $(0,\sqrt{E_{K^+}})$, by the Implicit Function Theorem.

We now make a first substitution $r=h(u)$ in the integral in $(\ref{period})$, obtaining
$$
T^+_{u_0}(E)=\frac{1}{\sqrt{2}}
\int_{\sqrt{E_{u_0}}}^{\sqrt{E}}\frac{dr}{h'(h^{-1}(r))\sqrt{E-r^2}}.
$$
A second substitution then presents itself naturally as $r=\sqrt{E}\sin \theta$, yielding
$$
T^+_{u_0}(E)=\frac{1}{\sqrt{2}}
\int_{\phi(E)}^{\pi/2} I(E,\theta) d \theta,
$$
where
$$
\phi(E):=\sin^{-1}(\sqrt{E_{u_0}}/\sqrt{E}),\textrm{ and } I(E,\theta):=\frac{1}{h'(h^{-1}(\sqrt{E}\sin \theta))}.
$$
Note for future reference that $\phi(E)\in (0,\pi/2)$, for all $E$ in $(E_{u_0},E_{K^+})$.

We now calculate the derivative of $T^+_{u_0}(E)$:
\begin{equation}\label{2-terms}
\frac{dT^+_{u_0}(E)}{dE}=\frac{1}{\sqrt{2}}
\left[\int_{\phi(E)}^{\pi/2} \frac{\partial I(E,\theta)}{\partial E} d \theta  - \frac{d\phi(E)}{dE}I(E,\phi(E)) \right].
\end{equation}
Differentiation under the integral sign is justified by the Leibniz integral rule, because for all $E$ in $(E_{u_0},E_{K^+})$ and all $\theta$ in 
$[\phi(E),\pi/2]$, 
$$
\left| \frac{\partial I (E,\theta)}{\partial E}\right| = \frac{1}{2\sqrt{E}} \left|   \frac{-h''(h^{-1}(\sqrt{E}\sin \theta))\sin \theta}{(h'(h^{-1}(\sqrt{E}\sin \theta)))^3}  \right| \leq  \frac{1}{2\sqrt{E}} \frac{M_E}{m_E^3},
$$
where $M_E$ is the positive maximum of the continuous function $|h''(u)|$ over the compact interval $[u_0,\beta(E)]$, and 
$m_E$ is the positive minimum of the continuous and positive function $h'(u)$ over the same interval.

It is easily verified that $d\phi(E)/dE<0$ and $I(E,\phi(E))>0$, for all $E$ in $(E_{u_0},E_K)$. Thus, to conclude the proof of this Lemma, it suffices to show that the integral appearing in the right hand side of $(\ref{2-terms})$ is non-negative. This integral is:
$$
\frac{1}{2\sqrt{E}}\int_{\phi(E)}^{\pi/2} - \frac{h''(h^{-1}(\sqrt{E}\sin \theta))\sin \theta}{(h'(h^{-1}(\sqrt{E}\sin \theta)))^3} d\theta,
$$
and we intend to evaluate it using an appropriate integration by parts. First, for any fixed $E$ in $(E_{u_0},E_{K^+})$, we set 
$$u(\theta)=h^{-1}(\sqrt{E}\sin \theta).
$$ 
Then the integrand of the latter integral factors as
$$
\frac{h''(u(\theta))}{(h'(u(\theta)))^3} \left(-\sin \theta d \theta \right) = \frac{h''(u(\theta))}{(h'(u(\theta)))^3} d(\cos \theta),
$$
and since
$$
\frac{du(\theta)}{d \theta}=\frac{\sqrt{E}\cos \theta}{h'(h^{-1}(\sqrt{E}\sin \theta))},
$$
integration by parts yields that
\begin{eqnarray}
\frac{1}{2\sqrt{E}}\int_{\phi(E)}^{\pi/2} - \frac{h''(u(\theta))\sin \theta}{(h'(u(\theta)))^3} d\theta = 
\frac{1}{2\sqrt{E}}\frac{h''(u(\theta))}{(h'(u(\theta)))^3} \cos \theta \Bigg|_{\theta=\phi(E)}^{\theta=\pi/2} + \nonumber \\
\frac{1}{2\sqrt{E}}\int_{\phi(E)}^{\pi/2}\cos \theta \frac{3(h''(u(\theta)))^2 - h'(u(\theta))h'''(u(\theta))}{(h'(u(\theta)))^4}  \frac{du(\theta)}{d\theta} d\theta \nonumber \\
= -\frac{1}{2\sqrt{E}}\frac{h''(u(\phi(E)))}{(h'(u(\phi(E))))^3} \cos (\phi(E)) + \label{bound-term} \\
\frac{1}{2}\int_{\phi(E)}^{\pi/2}\frac{3(h''(u(\theta)))^2 - h'(u(\theta))h'''(u(\theta))}{(h'(u(\theta)))^5} \cos^2 \theta d\theta. \label{int-term}
\end{eqnarray}
We claim that if ${\bf C1}^+$ and ${\bf C2}^+$ hold, then both quantities in $(\ref{bound-term})$ and in $(\ref{int-term})$ are non-negative. In order to prove this claim, we first derive two useful identities. These are obtained by repeatedly differentiating the identity $h^2=F^+$, see $(\ref{h-function})$, and expressing the derivatives for $h',h''$ and $h'''$ in terms $F^+,(F^+)',(F^+)''$ and $(F^+)'''$. One useful identity is
$$
h''=\frac{2F^+(F^+)'' - ((F^+)')^2}{4(F^+)^{3/2}}=(\sqrt{F^+})'',
$$
where the second equality in the previous line is easily established by direct calculation of the second derivative of $\sqrt{F^+}$.
This identity implies that if condition ${\bf C1}^+$ holds, then the boundary term in $(\ref{bound-term})$ is non-negative (recall that $h'>0$ and $\phi(E)\in (0,\pi/2)$).

A second useful identity is 
\begin{eqnarray*}
3(h'')^2-h'h''' &=& \frac{6F^+((F^+)'')^2 - 3((F^+)')^2(F^+)''-2F^+(F^+)'(F^+)'''}{8(F^+)^2}\\
&=&\frac{((F^+)')^4}{8(F^+)^2} \left(\frac{F^+}{((F^+)')^2} \right)'',
\end{eqnarray*}
where again, the last equality in the previous line is easily established by direct calculation of the second derivative of $F^+/((F^+)')^2$. 
Notice that the left hand side of this identity also appears in the 
numerator of the fraction found in the integrand of the integral in $(\ref{int-term})$. 
Consequently, if condition ${\bf C2}^+$ holds, then the integral term in $(\ref{int-term})$ is non-negative because its integrand is non-negative (recall again that $h'>0$). This proves our claim, and concludes the proof of the Lemma.
\end{proof}

\begin{figure}
    \centering
    \includegraphics[scale = 0.2]{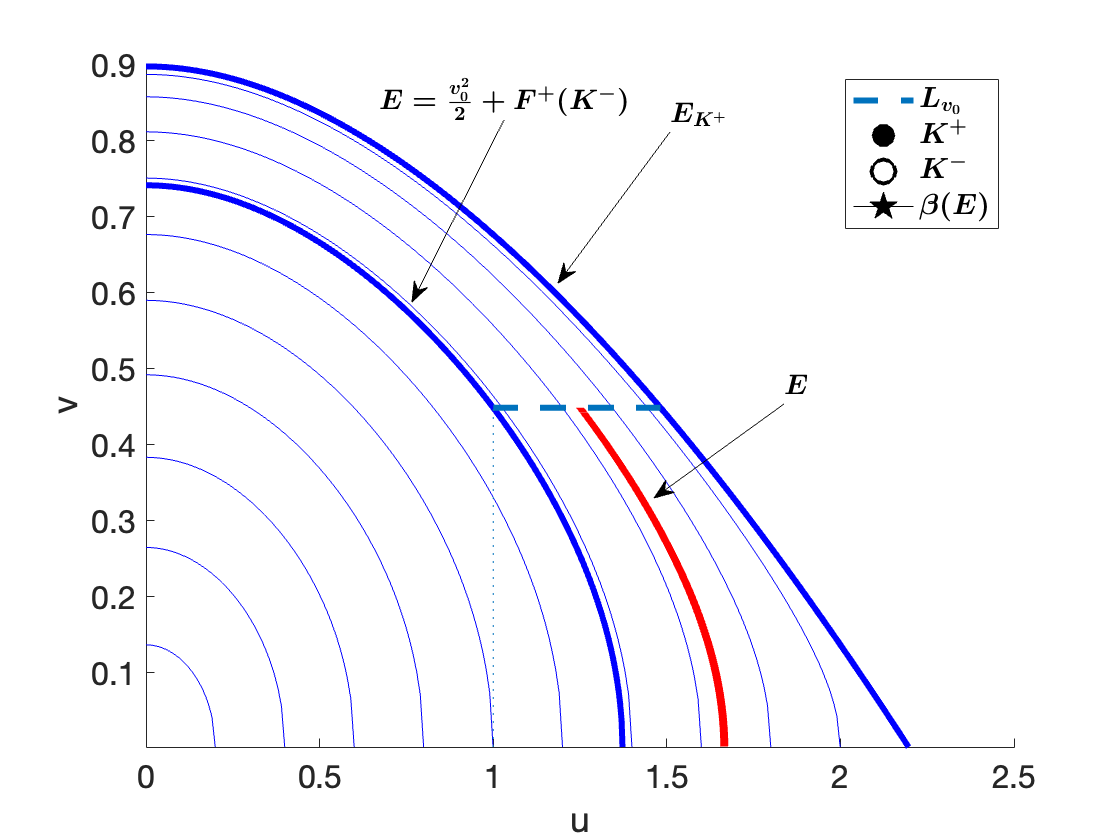}
    \caption{Orbits of the Hamiltonian system $(\ref{ham1})-(\ref{ham2})$. Reaction rate $f^+(u)=r^+u(1-u/K^+)$, with $r^+=1$, $K^+=2.2$ and diffusion constant $d^+=2$. Line segment $L_{v_0}$ for 
    $v_0=0.4491$.}
    \label{transversal_v0}
\end{figure}

Next, we again consider the Hamiltonian $H^+$, and we fix any $v_0$ in $(0,\sqrt{2E_{K^+}})$. 
This time we shall consider a family of solutions of $(\ref{ham1})-(\ref{ham2})$ whose initial conditions belong to the line segment 
$L_{v_0}=\{(u,v_0)\,|\, K^-<u< (F^+)^{-1}(E_{K^+}-v_0^2/2)\}$, as illustrated in Figure $\ref{transversal_v0}$. 
The vector field of $(\ref{ham1})-(\ref{ham2})$ is transversal to $L_{v_0}$. 
We parameterize the solutions starting on $L_{v_0}$ by the (constant) value $E$ the Hamiltonian $H^+$ takes on them. Since $F^+$ is strictly increasing over $(0,K^+)$, $E$ ranges over $(v_0^2/2 +F^+(K^-) , E_{K^+})$ as the initial condition of the solution varies in $L_{v_0}$. 
For each $E$ in $(v_0^2/2+F^+(K^-),E_{K^+})$, the corresponding solution intersects the $u$-axis in the point $(\beta(E),0)$, where, as before, $E\to \beta(E)$ is the unique $C^3$ function satisfying that $F^+(\beta(E))=E$.
For each $E$ in $(v_0^2/2+F^+(K^-),E_{K^+})$, let $T^+_{v_0}(E)$ 
denote the time it takes for the solution to reach the point $(\beta(E),0)$. 
%The map $T(E)$ is (at least) $C^1$ on the interval $(E_{u_0},E_K)$ because the vector field of system $(\ref{ham1})-(\ref{ham2})$ is $C^2$ in ${\cal R}$. In fact, 
Then $T^+_{v_0}(E)$ can be represented by the following Lebesgue integral:
\begin{equation*}
T^+_{v_0}(E)=\int_{u_0(E)}^{\beta(E)} \frac{du}{\sqrt{2(E-F^+(u))}},\textrm{ for all } E\in (v_0^2/2+F^+(K^-),E_{K^+}),
\end{equation*}
where
$$
u_0(E):=(F^+)^{-1}\left(E-v_0^2/2 \right),\textrm{ for all } E\in (v_0^2/2+F^+(K^-),E_{K^+}).
$$
The latter function $u_0(E)$ is a well-defined $C^3$ function because $F^+$ is strictly increasing on $(0,K^+)$. 
Using a similar argument as the one given before for the map $T_{u_0}^+(E)$, one can show that $T^+_{v_0}(E)<+\infty$ for all $E$ in $(v_0^2/2+F^+(K^-),E_{K^+})$.

We omit the proof of the following Lemma because it is entirely analogous to the proof of Lemma $\ref{monotone-period1}$.
\begin{lemma} \label{monotone-period2}
Suppose that ${\bf SA}$, ${\bf C1}^+$ and ${\bf C2}^+$ hold. Fix any $v_0$ in $(0,\sqrt{2E_{K^+}})$. 

Then $\frac{dT^+_{v_0}(E)}{dE}>0$, for all $E$ in $\left(v_0^2/2+F^+(K^-),E_{K^+} \right)$.
\end{lemma}

To prove our second main result Theorem $\ref{main2}$ (in which, when compared to Theorem $\ref{main}$, ${\bf M}^-$ is replaced by the pair of conditions ${\bf C1}^-$ and ${\bf C2}^-$), we must study certain solutions of the following Hamiltonian system with Hamiltonian $H^-$:
\begin{eqnarray}
u'&=& v, \label{ham3}\\
v'&=& -\frac{1}{d^-}f^-(u), \label{ham4}
\end{eqnarray}
defined in the region ${\cal R}=\{(u,v)\in \reals^2\,| \, u\geq 0\}$.

\begin{figure}
    \centering
    \includegraphics[scale = 0.2]{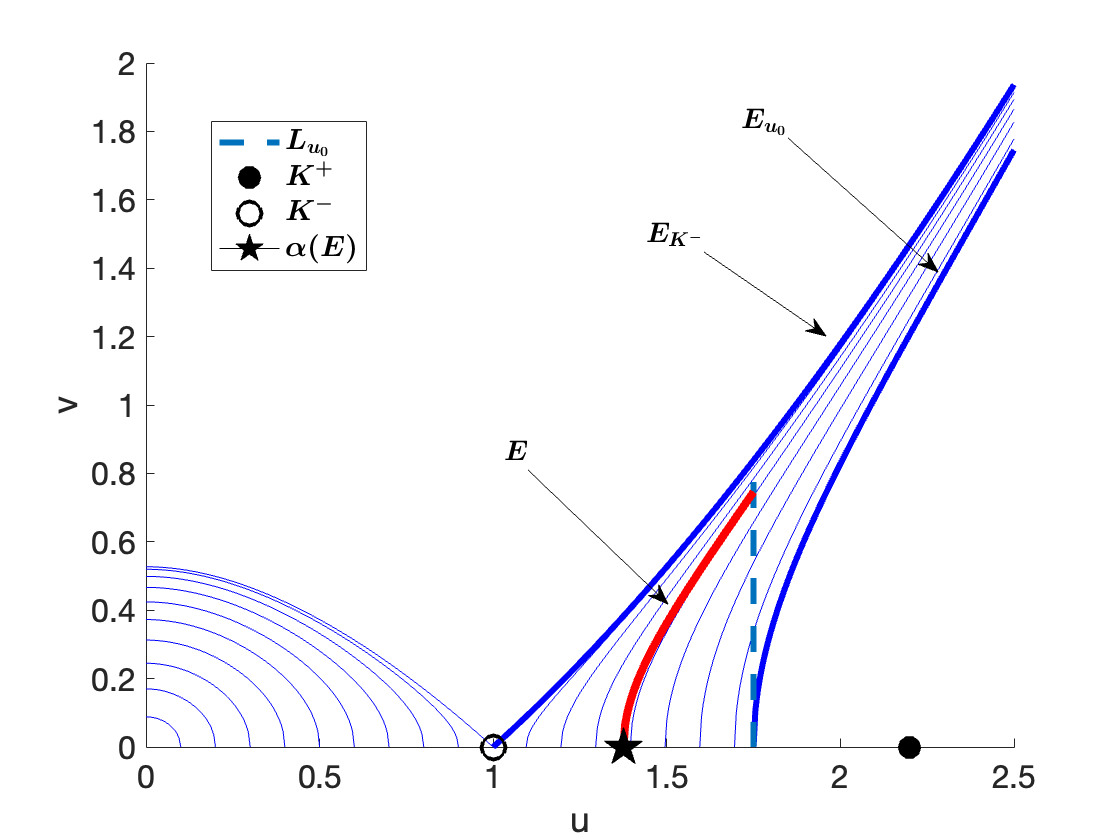}
    \caption{Orbits of the Hamiltonian system $(\ref{ham3})-(\ref{ham4})$. Reaction rate $f^-(u)=r^-u(1-u/K^-)$, with $r^-=1$, $K^-=1$ and diffusion constant $d^-=1.2$. Line segment $L_{u_0}$ for 
    $u_0=1.75$.}
    \label{transversalLem8_u0}
\end{figure}

Fix any $u_0$ in $(K^-,K^+)$ and set $E_{K^-}:=F^-(K^-)$ and $E_{u_0}:=F^-(u_0)$, 
and note that $E_{u_0}<E_{K^-}$ since $F^-$ is strictly decreasing on $(K^-,K^+)$. Define
\begin{equation*}
T_{u_0}^-(E)=\int_{\alpha(E)}^{u_0} \frac{du}{\sqrt{2(E-F^-(u))}},\textrm{ for all }E\in (E_{u_0},E_{K^-}),
\end{equation*}
where $E\to \alpha(E)$ is the unique $C^3$ function taking values in and onto $(K^-,u_0)$ such that $F^-(\alpha(E))=E$ (uniqueness follows from the fact that $F^-$ is strictly decreasing on $(K^-,K^+)$, by the Implicit Function Theorem). Then $T_{u_0}^-(E)$ represents the time it takes for the solution of the system $(\ref{ham3})-(\ref{ham4})$ starting at $(\alpha(E),0)$ to intersect the line where $u=u_0$, 
see Figure $\ref{transversalLem8_u0}$.

\begin{lemma} \label{monotone-period3}
Suppose that ${\bf SA}$, ${\bf C1}^-$ and ${\bf C2}^-$ hold. Fix any $u_0$ in $(K^-,K^+)$. 

Then $\frac{dT^-_{u_0}(E)}{dE}>0$, for all $E$ in $(E_{u_0},E_{K^-})$.
\end{lemma}
\begin{proof}
Define
$$
h(u):=\left( G^-(u)\right)^{1/2},\textrm{ for all }u\in (K^-,K^+),
$$
where $G^-(u)$ was defined in $(\ref{new-potential})$.

Then $h'(u)<0$ for all $u\in (K^-,K^+)$, and therefore the inverse $h^{-1}$ of $h$ exists. Also note that $h$ is $C^3$ on $(K^-,K^+)$, and so is $h^{-1}$, by the Implicit Function Theorem.

Defining $E_{K^+}:=F^-(K^+)$, and using the substitution $r=h(u)$, we find that
$$
T_{u_0}^-(E)=\frac{1}{\sqrt{2}} \int_{\sqrt{E-E_{K^+}}}^{\sqrt{E_{u_0}-E_{K^+}}} \frac{dr}{h'(h^{-1}(r))\sqrt{(E-E_{K^+})-r^2}}.
$$
Using the substitution $r=\sqrt{E-E_{K^+}} \sin \theta$, we obtain
$$
T_{u_0}^-(E)=\frac{1}{\sqrt{2}} \int_{\frac{\pi}{2}}^{\phi(E)}I(E,\theta)d\theta,
$$
where
$$
I(E,\theta):= \frac{1}{h'\left(h^{-1}\left(\sqrt{E-E_{K^+}}\sin \theta \right) \right)}, \textrm{ and }\phi(E):=\sin^{-1} \left( \frac{\sqrt{E_{u_0}-E_{K^+}}}{\sqrt{E-E_{K^+}}}    \right).
$$
Note that $\phi(E)$ belongs to $(0,\pi/2)$, and is strictly decreasing with respect to $E$.

Taking the derivative with respect to $E$, we find that
\begin{eqnarray*}
\frac{dT^-_{u_0}(E)}{dE}&=&\frac{1}{2\sqrt{2(E-E_{K^+})}} \int_{\frac{\pi}{2}}^{\phi(E)}\frac{-h''(h^{-1}(\sqrt{E-E_{K^+}}\sin \theta))\sin \theta}{\left(h'(h^{-1}(\sqrt{E-E_{K^+}\sin \theta})) \right)^3} d\theta \\
&+&\frac{1}{\sqrt{2}} \frac{d\phi(E)}{dE} I(E,\phi(E)).
\end{eqnarray*}
Since $h'<0$, the second term on the right hand side is positive, hence, to conclude the proof of the Lemma, it suffices to show that the integral is nonnegative. For each fixed $E$, we set $u(\theta)=h^{-1}(\sqrt{E-E_{K^+}}\sin \theta)$, and then an integration by parts shows that
\begin{eqnarray*}
&&\int_{\frac{\pi}{2}}^{\phi(E)}\frac{-h''(h^{-1}(\sqrt{E-E_{K^+}}\sin \theta))\sin \theta}{\left(h'(h^{-1}(\sqrt{E-E_{K^+}\sin \theta})) \right)^3} d\theta
=\frac{h''(u(\phi(E)))}{(h'(u(\phi(E))))^3}\cos (\phi(E))\\
&+&\sqrt{E-E_{K^+}} \int_{\pi/2}^{\phi(E)} \frac{3(h''(u(\theta)))^2 - h'(u(\theta))h'''(u(\theta))}{(h'(u(\theta)))^5} \cos^2 \theta d\theta.
\end{eqnarray*}
If ${\bf C1}^-$ holds, then the first term on the right hand side is non-negative (recall that $h'<0$ and that $\phi(E)\in (0,\pi/2)$). Similarly as in the proof of Lemma $\ref{monotone-period1}$, the numerator of the integrand of the integral on the right hand side is nonnegative when ${\bf C2}^-$ holds. Since $h'<0$ and $\phi(E)\in(0,\pi/2)$, it follows that the second term on the right hand side is nonnegative as well. This concludes the proof of this Lemma.
\end{proof}

\begin{figure}
    \centering
    \includegraphics[scale = 0.2]{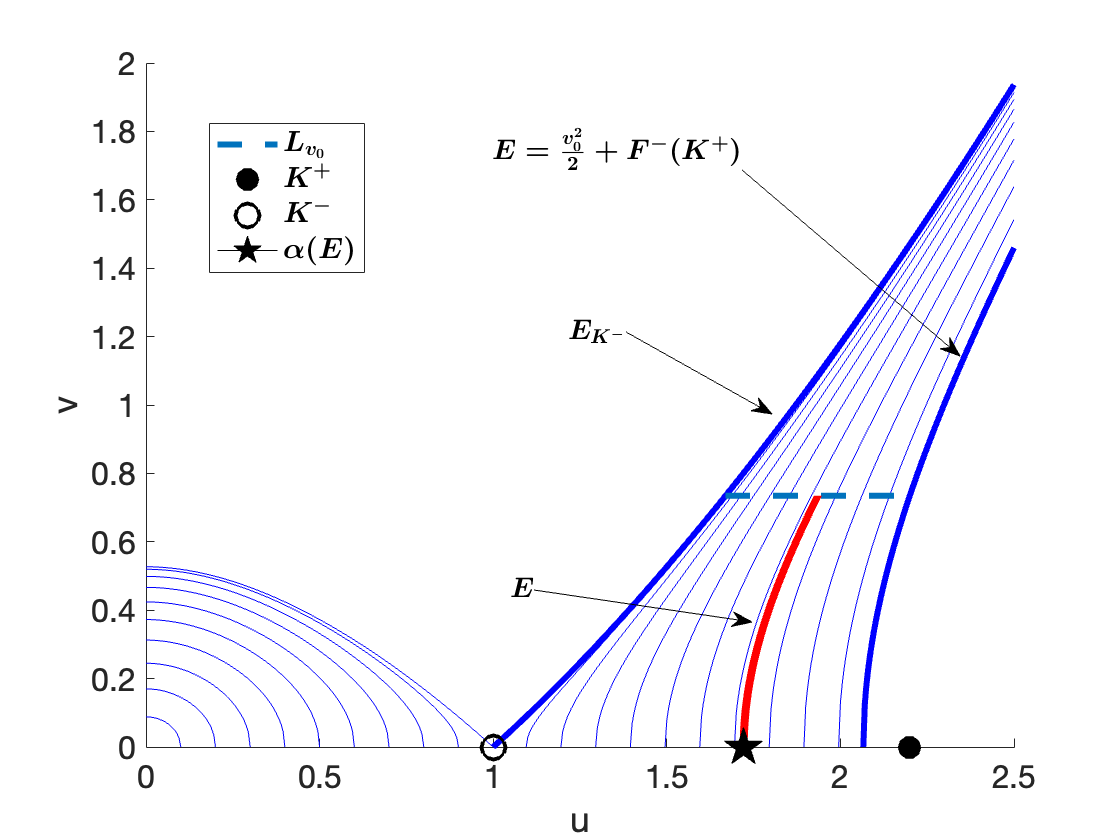}
    \caption{Orbits of the Hamiltonian system $(\ref{ham3})-(\ref{ham4})$. Reaction rate $f^-(u)=r^-u(1-u/K^-)$, with $r^-=1$, $K^-=1$ and diffusion constant $d^-=1.2$. Line segment $L_{v_0}$ for 
    $v_0=0.7348$.}
    \label{transversalLem9_v0}
\end{figure}

Let $E_{K^-}:=F^-(K^-)$, and $F^-_{\infty}:= \lim_{u\to +\infty} F^-(u)$. 
The latter limit exists (possibly as $-\infty$) because $F^-$ is decreasing on $(K^-,+\infty)$. 
Fix any $v_0$ in $\left(0,\sqrt{2(E_{K^-}-F^-_\infty)}\right)$. Define
\begin{equation*}
T_{v_0}^-(E)=\int_{\alpha(E)}^{u_0(E)} \frac{du}{\sqrt{2(E-F^-(u))}},\textrm{ for all }E\in \left(v_0^2/2+F^-(K^+),E_{K^-} \right), 
\end{equation*}
where, as before, $E\to \alpha(E)$ is the unique $C^3$ function such that $F^-(\alpha(E))=E$, and where 
$$
u_0(E):=(F^-)^{-1}\left(E-v_0^2/2 \right), \textrm{ for all }E\in \left(v_0^2/2+F^-(K^+),E_{K^-} \right).
$$
The latter function is a well-defined $C^3$ function since $F^-$ is strictly decreasing on $(K^-,+\infty)$. The function $T_{v_0}^-(E)$ represents the time it takes for the solution of the system $(\ref{ham3})-(\ref{ham4})$ starting at $(\alpha(E),0)$ to intersect the line where $v=v_0$, see Figure $\ref{transversalLem9_v0}$.

The proof of the following result is entirely analogous to that of Lemma $\ref{monotone-period3}$ and therefore omitted.
\begin{lemma} \label{monotone-period4}
Suppose that ${\bf SA}$, ${\bf C1}^-$ and ${\bf C2}^-$ hold. Fix any $v_0$ in 
\newline $\left(v_0^2/2+F^-(K^+),E_{K^-} \right)$. 

Then $\frac{dT^-_{v_0}(E)}{dE}>0$, for all $E$ in $\left(v_0^2/2+F^-(K^+),E_{K^-} \right)$.
\end{lemma}

\end{document}